%% file: optimal_linearization_arvix_v1.tex
\documentclass[journal,onecolumn]{IEEEtran}
\ifCLASSINFOpdf
  \usepackage[pdftex]{graphicx}
  \graphicspath{{./Figures/}}
\else
\fi
\usepackage{amsmath}
\usepackage{amsfonts}
\usepackage{amsthm}

\usepackage{multirow}
\usepackage{booktabs}
\usepackage[skip=0pt]{caption}

\newcommand{\setK}{$\mathbf{K}$}

\newcommand{\Rn}{$\mathbb{R}^n$}

\begin{document}
%
\title{Optimal Linearizations of Power Systems with Uncertain Supply and Demand}
%
%
%

\author{Marc Hohmann,
        Joseph Warrington
        and~John~Lygeros
\thanks{This research project is financially supported by the Swiss Innovation Agency Innosuisse and by NanoTera.ch under the project HeatReserves, and is part of the Swiss Competence Center for Energy Research SCCER FEEB\&D.}
\thanks{M. Hohmann is with the Urban Energy Systems Group, Empa, Swiss Federal Laboratories for Materials Science and Technology, \"Uberlandstrasse 129, 8600 D\"ubendorf , Switzerland.
{\tt\small marc.hohmann@empa.ch}.}
\thanks{J. Warrington and J. Lygeros are with the Automatic Control Laboratory, ETH Zurich, Physikstrasse 3, 8092 Z\"urich, Switzerland.}}

%
%

\markboth{}
{Hohmann \MakeLowercase{\textit{et al.}}: Optimal Linearizations of Power Systems}
%



\maketitle

\begin{abstract}
Linearized models of power systems are often desirable to formulate tractable control and optimization problems that still reflect real-world physics adequately under various operating conditions. In this paper, we propose an approach that can make use of known data concerning the distribution of demand, and/or intermittent supply, to minimize expected model inconsistency with respect to the original non-linear model. The optimal linearization is obtained by approximating a generalized moment problem with a hierarchy of sparse semi-definite relaxations. The output is a linearization point that minimizes the expected absolute constraint violation with respect to the uncertain supply and demand. Numerical results for different power systems networks demonstrate the accuracy and scalability of our linearization method. 
\end{abstract}


%
\IEEEpeerreviewmaketitle

\section{Introduction}
%
%
%
%
Fast, repeated power flow computations are required in a number of applications where supply and demand are uncertain. These include transmission expansion planning, day-ahead contingency analysis, and real-time load flow. However, the nonlinear constraints exhibited by power systems generally render the associated decision problems, most notably the Optimal Power Flow (OPF) problem, non-convex \cite{Taylor2015a}, making it difficult to solve exactly for large instances. This difficulty can be overcome by 1) the ad-hoc method of applying local optimization techniques such as interior-point methods to compute a locally optimal solution, 2) finding a convex relaxation of the non-convex optimization problem, or 3) linearizing the system model.

A convex relaxation constructs a relaxed convex set of constraints around the feasible operating range of the power system and provides a lower bound to the value of the optimal power flow problem. Under some conditions, exact global solutions can be extracted. The work of \cite{low2014b} details conditions under which second-order or semi-definite cone relaxations result in exact global solutions. Recent work \cite{Ghaddar2014}, \cite{Molzahn2015} and \cite{Josz2015} has focused on moment relaxations that can provide exact solutions of OPF cases for which these more common relaxations fail. Moment relaxations approximate a generalized moment problem (GMP) \cite{Lasserre2014}, the optimal solution of which is equivalent to the global optimum of the OPF, using semi-definite programming (SDP). As the computational complexity of moment relaxations increases rapidly with system size, the authors of \cite{Ghaddar2014} and \cite{Molzahn2015} applied sparsity exploitation techniques developed in \cite{Waki2006} to solve medium-sized OPF problems.

Linearization schemes are popular, because they lead to convex, albeit inexact, programs that scale favourably and can be solved very efficiently. They are typically preferred for contingency or probabilistic load flow studies that require a great number of analyses. Quantifying the approximation introduced by linearization is a challenging task that has only recently gained interest. In \cite{Molzahn2016}, the authors evaluate the worst case error of the well-known DC-flow approximation \cite{Christie2000}. The \emph{flat}-voltage approximation, a common starting point to derive application-specific linearized power flow equations such as the DC-flow, is shown to maintain active-power balance, but only under the restrictive assumption of a lossless network \cite{dhople2015}. The authors of \cite{Bolognani2016} and \cite{dhople2015} offer a voltage approximation linear in active and reactive demands, referred to as the \emph{No-Load} profile, for which bounds on the approximation error can be obtained as a function of the grid parameters.

If supply and demand are uncertain, it is desirable to linearize around a point that provides a good approximation under a range of conditions. Ideally, current statistical information should be used to adapt the linearization point to facilitate a range of online computations. To the best of our knowledge, no methods exploit statistical information about uncertain demand\footnote{For readability, we refer only to uncertain demand, but our approach is equally applicable to non-dispatchable sources such as wind and solar power.} to determine linearization points. Here we do this. By defining discrepancy as the summed expected magnitudes of constraint violations of the linearized problem with respect to a given probability distribution of electric loads, we provide a method for determining linearization points that minimize this quantity based on moment relaxations.

We formulate the full linearization problem as a two-stage stochastic polynomial optimization problem \cite{stoch2003}. In the first stage, a linearization point $x_0$ is selected such that the expected constraint violation at the second stage optimum is minimized. The second stage comprises an approximate, convex problem, in which nonlinear equality constraints have been linearized around $x_0$. Unfortunately, this two-stage setting is generally intractable because the optimality conditions, including complementary slackness conditions, and dual variables of the second stage have to be represented explicitly. 

Rather than solving this two-stage problem, we approximate the linearization point using SDP relaxations of a GMP inspired by \cite{Lasserre2010}. Specifically, we use the first moment, an output of the SDP relaxation, as a linearization point, and show that it minimizes the model discrepancy defined above. In addition, we show that in practice even low order SDP relaxations yield a satisfactory linearization point.

Standard moment methods \cite{Ghaddar2014}, \cite{Molzahn2015} and \cite{Josz2015} were developed to solve OPF problems \emph{online}, but their high computational cost is a problem as the computation has to be repeated many times. In contrast, the method we propose shifts the computational effort associated with the moment SDP \emph{offline}, where time requirements are less strict. The solution, namely a linearization point, is then used \emph{online} to formulate linearized OPF problems, which can be solved repeatedly and efficiently. We demonstrate that, thanks to the incorporation of demand statistics, the expected constraint violations are smaller in magnitude than under the method developed in \cite{dhople2015} and \cite{Bolognani2016}.


The OPF problem is presented in Section \ref{modelling}. We establish the full optimal linearization of the OPF problem in Section \ref{optimallinear}. Section \ref{approxoptimpoint} justifies the selection of the first moment of the parametric OPF problem as a linearization point and provides the dense and sparse SDP relaxations to approximate the first moment. The method is tested on standard benchmark networks in Section \ref{numericals}. We provide some concluding remarks in Section \ref{conclusion}.

\section{OPF Problem Formulation}\label{modelling}
Consider a power system described by a set of buses $\mathcal{N}=\{1,\ldots,N_B\}$ and a set $\mathcal{B}$ containing all pairs of buses connected by $N_L$ lines. The set $\mathcal{N}_G\subseteq\mathcal{N}$ denotes the $N_G$ generator buses. First, we state the set of power flow equations in rectangular voltage coordinates to obtain a polynomial description of the system:
\begin{subequations}
	\begin{align}
		 p_i(P_i,E,F):=&P_i-P_{Li}-\Bigg(\sum_{k=1}^{N_B} G_{ik}(E_i E_k+F_i F_k)\nonumber\\
		 &+B_{ik}(F_i E_k - E_i F_k)\Bigg)\nonumber\\
		 =&0,\quad\forall i\in\mathcal{N},\label{eq:nodeactive}\\
		q_i(Q_i,E,F):=&Q_i-Q_{Li}-\Bigg(\sum_{k=1}^{N_B} G_{ik}(F_i E_k - E_i F_k)\nonumber\\
		&-B_{ik}(E_i E_k + F_i F_k)\Bigg)\nonumber\\
		=&0,\quad\forall i\in\mathcal{N}\label{eq:nodereactive},
	\end{align}
\end{subequations}
where $E_i$ and $F_i$ are the real and imaginary components of the voltage phasor vectors $E\in\mathbb{R}^{N_B}$ and $F\in\mathbb{R}^{N_B}$, $P_i$ and $Q_i$ are the active and reactive components of the power generator injection vectors $P\in\mathbb{R}^{N_G}$ and $Q\in\mathbb{R}^{N_G}$, $P_{Li}$ and $Q_{Li}$ are the active and reactive components of the load vectors $P_{L}\in\mathbb{R}^{N_B}$ and $Q_{L}\in\mathbb{R}^{N_B}$; and $G_{ik}$ and $B_{ik}$ are the conductance and susceptance components of the admittance matrix $G+jB\in\mathbb{C}^{N_B\times N_B}$. The active and reactive power of the line power flows between two buses $(l,m)\in\mathcal{B}$ are given by:
\begin{subequations}
	\allowdisplaybreaks
	\begin{align}
	&P_{lm}=\frac{1}{\tau^2_{lm}}(g_{lm}+\frac{g_{sh,lm}}{2})(E_l^2+F_l^2)\nonumber\\
	&+\frac{1}{\tau_{lm}}(b_{lm}\sin(\theta_{lm})-g_{lm}\cos(\theta_{lm}))(E_l E_m+F_l F_m)\nonumber\\
	&+\frac{1}{\tau_{lm}}(g_{lm}\sin(\theta_{lm})+b_{lm}\cos(\theta_{lm}))(E_l F_m-F_l E_m),\label{eq:lineactive1}\\
	&P_{ml}=(g_{lm}+\frac{g_{sh,lm}}{2})(E_m^2+F_m^2)\nonumber\\
	&-\frac{1}{\tau_{lm}}(g_{lm}\cos(\theta_{lm})+b_{lm}\sin(\theta_{lm}))(E_l E_m+F_l F_m)\nonumber\\
	&+\frac{1}{\tau_{lm}}(g_{lm}\sin(\theta_{lm})-b_{lm}\cos(\theta_{lm}))(E_l F_m-F_l E_m),\label{eq:lineactive2}\\
	&Q_{lm}=-\frac{1}{\tau^2_{lm}}(b_{lm}+\frac{b_{sh,lm}}{2})(E_l^2+F_l^2)\nonumber\\
	&+\frac{1}{\tau_{lm}}(b_{lm}\cos(\theta_{lm})+g_{lm}\sin(\theta_{lm}))(E_l E_m+F_l F_m)\nonumber\\
	&+\frac{1}{\tau_{lm}}(g_{lm}\cos(\theta_{lm})-b_{lm}\sin(\theta_{lm}))(E_l F_m-F_l E_m),\label{eq:linereactive1}\\
	&Q_{ml}=-(b_{lm}+\frac{b_{sh,lm}}{2})(E_m^2+F_m^2)\nonumber\\
	&+\frac{1}{\tau_{lm}}(b_{lm}\cos(\theta_{lm})-g_{lm}\sin(\theta_{lm}))(E_l E_m+F_l F_m)\nonumber\\
	&+\frac{1}{\tau_{lm}}(g_{lm}\cos(\theta_{lm})+b_{lm}\sin(\theta_{lm}))(F_l E_m-E_l F_m),\label{eq:linereactive2}
	\end{align}
\end{subequations}
where $g_{lm}+jb_{lm}$ and $g_{sh,lm}+jb_{sh,lm}$ represent the series and total shunt admittance of the $\Pi$-model for the line from bus $l$ to $m$. Transformers, in series with the $\Pi$-model, are represented by the fixed complex turns ratio $1:\tau_{lm}e^{j\theta_{lm}}$ \cite{Zimmerman2011}. If there is no transformer, then $\tau_{lm}=1$ and $\theta_{lm}=0$. The line power flows $(P_{lm},P_{ml},Q_{ml},P_{ml})$ are denoted by $S$.
We consider the following OPF formulation: 
\begin{subequations}\label{eq:opf}
	\begin{align}
		\min_{E,F,P,Q,S} \quad &\sum_{k\in\mathcal{G}} f_k(P_k,Q_k)\\
	    \text{s.t.} \quad&\underline{P}_i \leq P_i \leq \overline{P}_i, \quad\forall i\in\mathcal{N}_G,\label{eq:activelimit}\\
	    &\underline{Q}_i \leq Q_i \leq \overline{Q}_i,\quad\forall i\in\mathcal{N}_G,\label{eq:reactivelimit}\\
	     &P_{lm}^2+Q_{lm}^2 \leq \overline{S}_{lm}^2, \quad\forall (l,m)\in\mathcal{B},\label{eq:lineab}\\
	    &P_{ml}^2+Q_{ml}^2 \leq \overline{S}_{lm}^2, \quad\forall (l,m)\in\mathcal{B},\label{eq:lineba}\\
	    & E_i^2+F_i^2 \leq \overline{V}_i^2, \quad\forall i\in\mathcal{N},\label{eq:voltageup}\\
	     &E_i^2+F_i^2 \geq \underline{V}_i^2, \quad\forall i\in\mathcal{N},\label{eq:voltagedown}\\
	     &E_1=1, F_1=0\label{eq:slack}\\
	     &\text{and (\ref{eq:nodeactive}), (\ref{eq:nodereactive}), (\ref{eq:lineactive1}), (\ref{eq:lineactive2}),  (\ref{eq:linereactive1}), (\ref{eq:linereactive2})},\nonumber
	\end{align}
\end{subequations}
where the functions $f_k(P_k,Q_k)$ are assumed to be convex. The loads $P_L$ and $Q_L$ are uncertain disturbances. The OPF (\ref{eq:opf}) is assumed to be solved with recourse, i.e.~the decision $(E,F,P,Q,S)$ are taken after the realization of $P_L$ and $Q_L$. Thus, the OPF solutions are parametric in $P_L$ and $Q_L$ that appear as affine terms in (\ref{eq:nodeactive}) and (\ref{eq:nodereactive}). Constraints (\ref{eq:activelimit}), (\ref{eq:reactivelimit}), (\ref{eq:voltageup}) and (\ref{eq:voltagedown}) limit the active and reactive power generation and the voltage magnitudes at each bus. The apparent line power limits are enforced with constraints (\ref{eq:lineab}) and (\ref{eq:lineba}) and expressions (\ref{eq:lineactive1})-(\ref{eq:linereactive2}). Constraints (\ref{eq:slack}) define the slack bus (in p.u.). To obtain convex inequality and non-convex equality constraints, we replace the constraints (\ref{eq:voltagedown}) by 
\begin{subequations}
	\begin{align}
		&X_i=E_i^2+F_i^2,\quad\forall i\in\mathcal{N},\label{eq:voltmag}\\
		&X_i \geq \underline{V}_i^2,\quad\forall i\in\mathcal{N},\label{eq:voltdown2}
	\end{align}
\end{subequations} where $X_i$ is the squared voltage phasor magnitude at bus $i\in\mathcal{N}$. Note that our formulation (\ref{eq:opf}) contains some redundant variables but it has been chosen so that it only involves polynomial constraints of degree at most 2. The reason for this is made clear in Section \ref{approxoptimpoint} and the Appendix.

\section{Optimal linearization}\label{optimallinear}
In this section, we derive the full optimal linearization of (\ref{eq:opf}). To streamline the notation, we write (\ref{eq:opf}) as a polynomial optimization problem:
\begin{subequations}\label{eq:parametricoptim}
	\begin{align}
	\min_{x \in \mathbb{R}^n} &f(x) \quad\label{eq:energy systemsa}\\
	\text{s.t.} \quad &h_i(x,y) = 0, \quad i=1,\ldots,N_h,\label{eq:energy systemsb}\\ 
	\quad &g_j(x) \geq 0, \quad j=1,\ldots,N_g,\label{eq:energy systemsc}
	\end{align}
\end{subequations}
where the vector of control decisions $x\in\mathbb{R}^{n}$ represents the OPF decisions $(E,F,P,Q,S,X)$ and the vector of parameters $y\in\mathbb{R}^{p}$ represents the loads $(P_L,Q_L)$. The constraints (\ref{eq:energy systemsb}) represent all the equality constraints (\ref{eq:nodeactive}), (\ref{eq:nodereactive}), (\ref{eq:lineactive1})-(\ref{eq:linereactive2}) and (\ref{eq:voltmag}). Note that $y$ enters the equality constraints as affine terms. The inequality constraints (\ref{eq:activelimit})-(\ref{eq:voltageup}) and (\ref{eq:voltdown2}) are represented by (\ref{eq:energy systemsc}). There are $N_h=3N_B+4N_L$ equality and  $N_g=2N_G+2N_B+2N_L$ inequality constraints. All generator cost functions are summarized with a function $f: \mathbb{R}^n\rightarrow\mathbb{R}$. Owing to the nature of real electricity demand, we assume that $y$ is restricted to the compact semi-algebraic set $\mathbf{Y}:=\{y\in\mathbb{R}^p: g_j(y)\geq 0,j=N_g+1,\ldots,N_g'\}$. The set of all possible combinations of control decisions and consumer demand is defined as
\begin{equation}
	\begin{aligned}
	\mathbf{K}:=\{&(x,y)\in\mathbb{R}^{n}\times\mathbb{R}^p,y\in\mathbf{Y}:\\
	& h_i(x,y)=0, i=1,\ldots,N_h;\\
	&\,g_j(x)\geq 0, j=1,\ldots,N_g\}.
	\end{aligned}
\end{equation}
The technical bounds (\ref{eq:activelimit})-(\ref{eq:voltagedown}) ensure that $\mathbf{K}$ is compact.

With the OPF formulation of Section \ref{modelling}, the difficulty in solving (\ref{eq:parametricoptim}) for a given $y$ arises only from the non-convex equality constraints (\ref{eq:energy systemsb}). Here we aim to approximate (\ref{eq:parametricoptim}) by a convex optimization problem:
\begin{subequations}\label{eq:convexparam}
	\begin{align}
	\min_{x_{\rm lin} \in \mathbb{R}^n} &f(x_{\rm lin})\\
	\text{s.t.} \quad &h^{\rm lin}_i(x_{\rm lin},y) = 0, \quad i=1,\ldots,N_h,\label{eq:convexparamb}\\ 
	&g_j(x_{\rm lin}) \geq 0, \quad j=1,\ldots,N_g,\label{eq:convexparamc}
	\end{align}
\end{subequations}
where $h^{\rm lin}_i(x_{\rm lin},y)$ are linearizations of $h_i(x,y)$ around an operating point $x_0$ and defined as $h^{\rm lin}_i(x_{\rm lin},y)=h_i(x_0,y)+\nabla_{x} h_i(x_0)^T (x_{\rm lin}-x_0)$, where $\nabla_{x} h_i(x_0)$ is the gradient of $h_i$ at $x_0$ with respect to $x$. Since $h_i(x,y)$ are affine in $y$ for problem (\ref{eq:opf}), only a linearization around $x$ must be determined.

For each $i=1,\ldots,N_h$, we define the signed linearization error as:
\begin{equation}\label{eq:signedlinearization}
	\epsilon_i(x_{\rm lin},y):=h_i(x_{\rm lin},y)-h^{\rm lin}_i(x_{\rm lin},y).
\end{equation}
 Let $x^*_{\rm lin}(y) \in \mathbb{R}^n$ denote an optimal parametric solution to (\ref{eq:convexparam}) for a given $y$. Expressing the Karush-Kuhn-Tucker (KKT) conditions of (\ref{eq:convexparam}) in terms of $h_i(x^*_{\rm lin},y)$ and $\epsilon_i(x^*_{\rm lin},y)$, we obtain:
\begin{equation}\label{eq:kkt}
\begin{aligned}
&\sum_{i=1}^{N_h} \lambda_i(y) \nabla_x \epsilon_i(x^*_{\rm lin}(y),y)=\\
&\quad\nabla_x f(x^*_{\rm lin}(y)) +\sum_{j=1}^{N_g} \kappa_j(y) \nabla_x g_j(x^*_{\rm lin}(y))\\
&\quad+\sum_{i=1}^{N_h} \lambda_i(y) \nabla_x h_i(x^*_{\rm lin}(y)),\\
&h_i(x^*_{\rm lin}(y),y) = \epsilon_i(x^*_{\rm lin}(y),y), \quad i=1,\ldots,N_h,\\ 
&g_j(x^*_{\rm lin}(y)) \geq 0,\quad\kappa_j(y) \geq 0,\quad j=1,\ldots,N_g,\\ 
&g_j(x^*_{\rm lin}(y))\kappa_j(y)=0, \quad j=1,\ldots,N_g,
\end{aligned}
\end{equation}
where $\lambda(y) \in \mathbb{R}^{N_h}$ and $\kappa(y) \in \mathbb{R}^{N_g}$ are the Lagrange multipliers of (\ref{eq:convexparamb}) and (\ref{eq:convexparamc}).

By minimizing $\sum_i |\epsilon_i(x^*_{\rm lin}(y),y)|$, the equality constraint violations can be minimized for an approximating solution $x^*_{\rm lin}(y)$.
Given a probability measure $\varphi$ on $\mathbf{Y}$, describing the uncertainty of the loads $(P_L,Q_L)$, we define the expected constraint violation as:
\begin{equation}\label{eq:expectedconviol}
\mathbf{E}_\varphi\Bigg(\sum_{i=1}^{N_h} |\epsilon_i(x_{\rm lin}^*(y),y)|\Bigg)=\int_{\mathbf{Y}} \sum_{i=1}^{N_h}|\epsilon_i(x_{\rm lin}^*(y),y)|d\varphi.
\end{equation} 

We aim to select the operating point $x_0$ such that the solutions of (\ref{eq:convexparam}) minimize the expected constraint violations. For this purpose, the optimal selection of a linearization point $x_0$ is cast as a two-stage stochastic problem (with respect to the distribution $\varphi$):
\begin{equation}
\begin{gathered}\label{eq:fullLinearization}
\min_{x_0 \in \mathbb{R}^n} \mathbf{E}_\varphi(v(x_0,y)),\\
\end{gathered}
\end{equation}
where $v(x_0,y)$ is the optimal value function of the second stage problem:
\begin{equation}
\begin{gathered}
v(x_0,y)=\min_{x_{\rm lin} \in \mathbb{R}^n,\lambda \in \mathbb{R}^{N_h},\kappa \in \mathbb{R}^{N_g}} \sum_{i=1}^{N_h} |\epsilon_i(x_{\rm lin},y)| \\
\quad \text{s.t. (\ref{eq:kkt})}.
\end{gathered}
\end{equation}
Note that the optimality constraints (\ref{eq:kkt}) in problem (12) encode the fact that we wish to penalize constraint violations at an economic optimum (according to the convexified model) after the linearization around $x_0$ has been carried out. However, problem (\ref{eq:fullLinearization}) is intractable even for small instances due to the additional variables and polynomial constraints of the KKT conditions, including complementarity constraints. To address this difficulty, we propose an approximate linearization procedure.
\section{Approximation of the optimal linearization point}\label{approxoptimpoint}
Assuming that the linearized optimal solutions differ only moderately from the global OPF solutions, we make an approximation by minimizing the linearization error with respect to the distribution of global OPF solutions:
\begin{equation}
\begin{aligned}\label{eq:meanconstraint}
\min_{x_0} \int_{\mathbf{Y}}& \sum_{i=1}^{N_h}|\epsilon_i(x^*(y),y)| d\varphi\\
\end{aligned}
\end{equation}
where $x^*(y) \in \mathbb{R}^n$ denotes an optimal parametric solution to (\ref{eq:parametricoptim}) for a given $y$. We show in the Appendix that using the expected value of the OPF solutions $\mathbf{E}_{\varphi}(x^*(y))$ as a linearization point is equivalent to minimizing (\ref{eq:meanconstraint}) in the case of constraints (\ref{eq:voltmag}) or minimizing a convex upper bound to (\ref{eq:meanconstraint}) in the case of constraints (\ref{eq:nodeactive}), (\ref{eq:nodereactive}) and (\ref{eq:lineactive1})-(\ref{eq:linereactive2}). Motivated by this, we approximate the expected value of OPF solutions $\mathbf{E}_{\varphi}(x^*(y))$ using SDP relaxations.
The characterization of optimal solutions of (\ref{eq:parametricoptim}) in the form of a GMP is presented in Section \ref{generalizedmomentproblem}, followed by its SDP relaxation in Section \ref{sdprelax}. Since the computational complexity of the SDP relaxation increases substantially with the system size, we discuss sparsity exploiting SDP relaxations in Section \ref{sparse}.
\subsection{Generalized Moment Problem}\label{generalizedmomentproblem}

The set $\mathcal{M}(\mathbf{K})_+$ denotes the non-negative Borel measures on $\mathbf{K}$. Let $\mu$ be a Borel probability representing the joint distribution of OPF decisions $x^*(y)$ and parameters $y$ denoted by $\mu\in\mathcal{M}(\mathbf{K})_+$, with $\mu(\mathbf{K})=1$. Let  $\pi: \mathcal{M}(\mathbf{K})_+ \rightarrow \mathcal{M}(\mathbf{Y})_+$ be the projection on the parameter set, defined by ($\pi\mu$)($\mathit{B}$)=$\mu$((\Rn $\times \mathit{B})~\cap$ \setK) for all Borel subsets $\mathit{B}$ of $\mathbf{Y}$. The following GMP based on \cite{Lasserre2010} encodes all instances of (\ref{eq:parametricoptim}) when the load uncertainty measure $\varphi$ is given:
\begin{equation}\label{GMBM}
\begin{gathered}
	\rho=\min_{\mu \in \mathcal{M}(\mathbf{K})_+}\int_{\mathbf{K}} fd\mu\quad\text{s.t.}\quad\pi\mu=\varphi,\\
\end{gathered}
\end{equation}
where the constraint $\pi\mu=\varphi$ imposes the load distribution $\varphi$ onto the parameters $y$. Problem (\ref{GMBM}) solves an optimization problem for all OPF instances of $(P_L,Q_L)$ while weighting the occurrence of $(P_L,Q_L)$ by $\varphi$. The optimal probability distribution $\mu^*$ of (\ref{GMBM}) characterizes the distribution of optimal parametric OPF solutions $(E^*,F^*,P^*,Q^*,S^*,X^*)$ when $(P_l,Q_l)$ follows $\varphi$. The explicit mapping $x^*(y)$ is not needed to obtain the expected value $\mathbf{E}_{\varphi}(x^*(y))$ because we have the following result by (\cite{Lasserre2010}, Corollary 2.3):
\begin{equation}\label{eq:parametricsol}
 \mathbf{E}_{\varphi}(x^*(y))=\int_{\mathbf{Y}} x^*(y) d\varphi=\mathbf{E}_{\mu^*}(x)=\int_{\mathbf{K}} x d\mu^*.
\end{equation}

Problem (\ref{GMBM}) is not tractable because it is infinite dimensional. For this reason, we formulate a hierarchy of SDP relaxations to approximate the moments of the distribution $\mu$. 
\subsection{SDP Relaxations of the GMP}\label{sdprelax}
As shown in \cite{Lasserre2010}, one can obtain a hierarchy of SDP relaxations approximating (\ref{GMBM}) by considering a finite number of moments of $\mu$. The relaxation degree $k\in\mathbb{N}$ provides a trade-off between the accuracy of the approximation and computational complexity. In the case of problem (\ref{eq:parametricoptim}), of which the constraints are at most second order polynomials, $k$ must satisfy $k\geq \lceil \textrm{deg}\,f/2 \rceil $, where ``$\textrm{deg}$'' is the degree of a polynomial and $\lceil b \rceil$ denotes the ceiling of a real number $b$ (smallest integer greater than or equal to the number $b$).


At level $k$ of the hierarchy, an SDP relaxation involving all the moments up to order $2k$ is solved. Let $\alpha=(\alpha_1,\ldots,\alpha_n)\in\mathbb{N}^n$ and $\gamma=(\gamma_1,\ldots,\gamma_p)\in\mathbb{N}^p$ be the integer vectors of dimension $n$ and $p$ respectively (to serve as multi-indices) and define $m_{\alpha\gamma}$ as the moments of the probability measure $\mu$ on $\mathbf{K}$ by:
\begin{equation}
m_{\alpha\gamma}:=\int_{\mathbf{K}} x^\alpha y^\gamma  d\mu,
\end{equation}
where the shorthand $x^\alpha=x_1^{\alpha_1} x_2^{\alpha_2}\ldots x_n^{\alpha_n}$ and $y^\gamma=y_1^{\gamma_1} y_2^{\gamma_2}\ldots y_p^{\gamma_p}$ are used for the monomials. Note that $m_{00}=1$ because $\mu$ is a probability measure. Let $\mathbf{m}_k$ be a vector containing all the moments $m_{\alpha\gamma}$ up to degree $2k$ such that $\sum_n \alpha_n+\sum_p \gamma_p \leq 2k$. The fixed moments of the known marginal probability measure $\varphi$ are denoted by $z_\gamma=\int_{\mathbf{Y}}y^\gamma d\varphi$, with $z_{0}=1$.
Let $u\in\mathbb{R}[x,y]$ be a polynomial, where $\mathbb{R}[x,y]$ is the ring of polynomials in $(x,y)$. We associate a linear mapping $L_{\mathbf{m}}:\mathbb{R}[x,y]\rightarrow\mathbb{R}$ with the moment vector $\mathbf{m}_k$ defined by:
\begin{equation}
	u=\sum_{\alpha\gamma}u_{\alpha\gamma}x^\alpha y^\gamma \mapsto L_{\mathbf{m}}(u)=\sum_{\alpha\gamma}u_{\alpha\gamma}m_{\alpha\gamma}.
\end{equation}
It can be shown that, if we restrict our attention to polynomials of degree $2k$, $L_\mathbf{m}$ defines a positive definite matrix $M_k(\mathbf{m})$, the so called moment matrix; if $u$ is a polynomial of degree $k$, applying $L_\mathbf{m}$ to $u^2$ leads to
\begin{equation}
L_{\mathbf{m}}(u^2)=\mathbf{u}^T M_k(\mathbf{m}) \mathbf{u},
\end{equation}
where $\mathbf{u}$ is the vector of coefficients of $u$ and $M_k(\mathbf{m})$ comprises entries $m_{\alpha\gamma}$ of $\mathbf{m}_k$. 
Since $u^2$ is non-negative it is easy to see that the moment matrix is symmetric positive semi-definite. The positive semi-definite localizing matrices $M_{k-1}(g_j\mathbf{m})$, used to enforce constraints $g_j\geq 0, j=1,\ldots,N_g'$, are derived in a similar fashion:
\begin{equation}
L_{\mathbf{m}}(g_ju^2)=\mathbf{u}^T M_{k-1}(g_j\mathbf{m}) \mathbf{u}\geq 0.
\end{equation}
The entries of $M_{k-1}(g_j\mathbf{m})$ are linear combinations of the moments $m_{\alpha\gamma}$. Finally, the localizing matrices of equality constraints $h_i=0, i=1,\ldots,N_h$ are defined as
\begin{equation}
L_{\mathbf{m}}(h_i u^2)=\mathbf{u}^T M_{k-1}(h_i\mathbf{m}) \mathbf{u}.
\end{equation}
Since an equality can be treated as two reverse inequalities, all entries of $M_{k-1}(h_i\mathbf{m})$, linear combinations of the moments $m_{\alpha\gamma}$, must be equal to zero. We write this as $M_{k-1}(h_i\mathbf{m})=0$.

Treating the equalities as two reverse inequalities, it is shown in \cite{Lasserre2014} that the sequence of moments $\mathbf{m}_k$ has a representing finite Borel measure $\mu$ on $\mathbf{K}$ if and only if the moment matrix and the localizing matrices are positive semi-definite for all $k\in\mathbb{N}$. Based on this, the relaxation of (\ref{GMBM}) of degree $k$ in the hierarchy involves solving an SDP of the form:
\begin{subequations}\label{eq:primalsdp}
	\begin{align}
	\rho_k=&\min_{\mathbf{m}_k}\enskip L_{\mathbf{m}}(f)\label{eq:primalsdpobj}\\
	\text{s.t.}\quad &\medspace M_k(\mathbf{m})\succeq 0\label{eq:primalsdppsd}\\
	&M_{k-1}(h_i\mathbf{m})=0, \quad i={1,\ldots,N_h},\\
	&M_{k-1}(g_j\mathbf{m})\succeq0, \quad j={1,\ldots,N_g'},\label{eq:primalsdplocal}\\
	&m_{0\gamma}=z_{\gamma}, \quad\forall\gamma\in \mathbb{N}_{2k}^{p},\label{eq:primalsdpfixed}
	\end{align}
\end{subequations}
where $\succeq 0$ requires a matrix to be positive semi-definite and $\mathbb{N}_{2k}^{p}$ denotes the set of all $\gamma$ such that $\sum_p \gamma_p\leq 2k$.
The objective (\ref{eq:primalsdpobj}) represents the expected cost of $f$ of (\ref{eq:energy systemsa}) with respect to a probability measure $\mu$ supported on $\mathbf{K}$. The constraint set (\ref{eq:primalsdpfixed}) imposes the known moments of $\varphi$ onto the marginal moments of $\mu$. 

The relaxation hierarchy provides increasingly accurate approximations of the moments of the distribution of optimal solution of (\ref{eq:primalsdp}) \cite{Lasserre2010}, i.e.~$\lim_{k\rightarrow\infty} \rho_k=\rho$, with $\rho_k\leq\rho_{k+1}$ for all $k$ and $\lim_{k\rightarrow\infty} \mathbf{m}_k=\mathbf{m}$, where $\mathbf{m}$ is the moment vector of $\mu^*$. Here we propose to use the first moment of $\mathbf{m}_k$ computed by (\ref{eq:primalsdp}), approximating the expected value $\mathbf{E}_\varphi(x^*(y))$, as a linearization point $x_0$ to obtain a convex approximation (\ref{eq:convexparam}).
\subsection{Sparsity exploiting SDP relaxation}\label{sparse}
Solving an SDP relaxation at level $k$ involves moment matrices with dimensions up to $\begin{pmatrix}
n+p+k\\k
\end{pmatrix}\times\begin{pmatrix}
n+p+k\\k
\end{pmatrix}$.
This means that, in a naive implementation, the matrices grow very quickly in dimension as a function of the number of buses and degree of relaxation. Since the dimension of the semi-definite constraints is typically the main bottleneck for SDP solvers, the authors of \cite{Waki2006} proposed SDP relaxations with more, but significantly smaller, moment and localizing matrices than (\ref{eq:primalsdp}) by exploiting the sparsity structure of the set $\mathbf{K}$ and the polynomial $f$.

Let $J$ be the set of all monomials contained in $h_1(x,y),\ldots,h_{N_h}(x,y)$, $g_1(x),\ldots,g_{N_g'}(y)$ and $f(x)$. A subset of monomials of $J$ with index $s\in\{1,\ldots,N_I\}$ only involves a subset of the variables $\{x_1,\ldots,x_n,y_1,\ldots,y_p\}$. Let $I_s$ define the subset of variables in question, and $n_s$ and $p_s$ the cardinality with respect to $x$ and $y$, respectively. Define the set $\mathcal{H}_{s} \subset \{1,\ldots,N_h\}$ as the set of indices of constraint functions $h_1(x,y),\ldots,h_{N_h}(x,y)$ that include at least one variable of $I_s$. Furthermore, define $\mathcal{G}_{s}\subset\{1,\ldots,N_g'\}$ as the set of indices of constraint functions $g_1(x),\ldots,g_{N_g'}(y)$ that include at least one variable of $I_s$.

Given a collection $\{I_1,\ldots,I_{N_I}\}$, a multi-measures moment problem can be formulated that, by virtue of (\cite{Lasserre2014}, Theorem 4.6), is equivalent to solving (\ref{GMBM}) if the collection $\{I_1,\ldots,I_{N_I}\}$ satisfies the \emph{running intersection property} defined by:

For each $s=1,\ldots,N_I-1$, we have
\begin{equation}\label{eq:RIP}
	I_{s+1}\cap\bigcup_{t=1}^s I_t\subseteq I_q \enskip\text{for some}\enskip q\leq s.
\end{equation}


Consider the following sparse SDP relaxation at level $k$ in the hierarchy:
\begin{subequations}\label{eq:sparsesdp}
	\begin{align}
	\upsilon_k =&\min_{\mathbf{m}}\enskip L_{\mathbf{m}}(f)\label{eq:sparsesdpobj}\\
	\text{s.t.}\quad & M_k(\mathbf{m},I_s)\succeq 0\label{eq:sparsesdppsd},\quad s=1,\ldots,N_I,\\
	&M_{k-1}(h_i\mathbf{m},I_s)=0,\\ &\quad\quad \forall i\in \mathcal{H}_s,s=1,\ldots,N_I,\\
	&M_{k-1}(g_j\mathbf{m},I_s)\succeq0,\\ &\quad\quad \forall j\in\mathcal{G}_s,s=1,\ldots,N_I,\label{eq:sparsesdplocal}\\
	&m_{0\gamma}(I_s)=z_{\gamma}(I_s), \quad\forall\gamma\in \mathbb{N}_{2k}^{p_s}, s=1,\ldots,N_I,\label{eq:sparsesdpfixed}
	\end{align}
\end{subequations}
where $m_{0\gamma}(I_s)$ are the moments of the parameters indexed by $I_s$, $M_k(\mathbf{m},I_s)$ are moment matrices constructed from the variables and parameters indexed by $I_s$; and $M_{k-1}(h_i\mathbf{m},I_s)$ and $M_{k-1}(g_j\mathbf{m},I_s)$ are localizing matrices for constraints indexed by $\mathcal{H}_s$ and $\mathcal{G}_s$ and constructed from the variables and parameters indexed by $I_s$. Contrary to (\ref{eq:primalsdpfixed}), we only impose the moments of $\varphi$ indexed by $I_s$, denoted as $z_{\gamma}(I_s)$, in (\ref{eq:sparsesdpfixed}).
By virtue of (\cite{Lasserre2014}, Theorem 4.7), if the collection $\{I_1,\ldots,I_{N_I}\}$ satisfies (\ref{eq:RIP}), the sparse SDP relaxation (\ref{eq:sparsesdp}) converges to the optimal solution of $(\ref{GMBM})$, i.e.~$\lim_{k\rightarrow\infty} \upsilon_k=\rho$.

An efficient method for identifying a collection $\{I_1,\ldots,I_{N_I}\}$ satisfying (\ref{eq:RIP}) was proposed in \cite{Waki2006} and used for power systems in \cite{Ghaddar2014}, \cite{Molzahn2015} and \cite{duan2017}. It is based on the chordal structure in a graph with vertices that correspond to variables and edges representing their interaction. There is an edge between to vertices if the variables appear in the same constraint or monomial of $f$.

Solving the sparse SDP relaxation (\ref{eq:sparsesdp}) at level $k$ involves matrices with dimensions $\begin{pmatrix}
n_s+p_s+k\\k
\end{pmatrix}\times\begin{pmatrix}
n_s+p_s+k\\k
\end{pmatrix}$. Thus, the smaller the cardinalities of all subsets $I_s$ found by the sparsity detection method, the lower the computational load of (\ref{eq:sparsesdp}).

\begin{figure}[!t]
	\centering
	\includegraphics[scale=1]{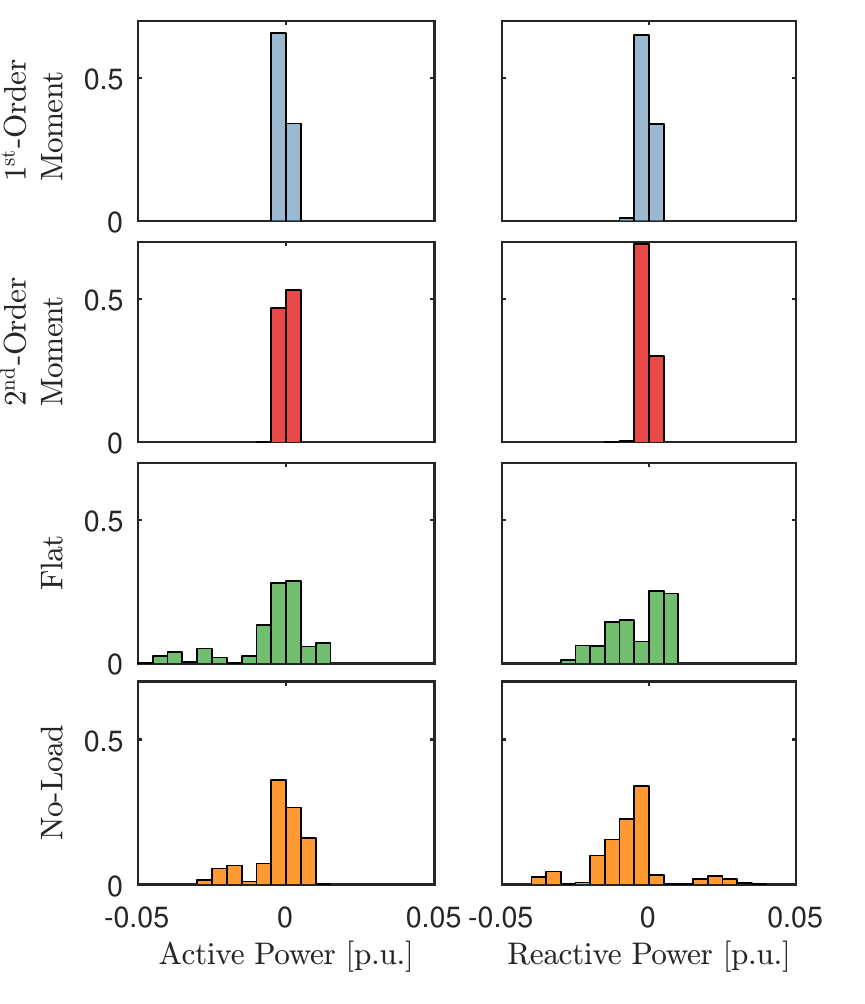}
	\caption{IEEE 14-Bus: Histograms of signed active and reactive power balance errors $p_i$ and $q_i$ for the four linearized models over all buses and 1000 scenarios.}\label{fig:case14} 
\end{figure}
\begin{figure}[!t]
	\centering
	\includegraphics[scale=1]{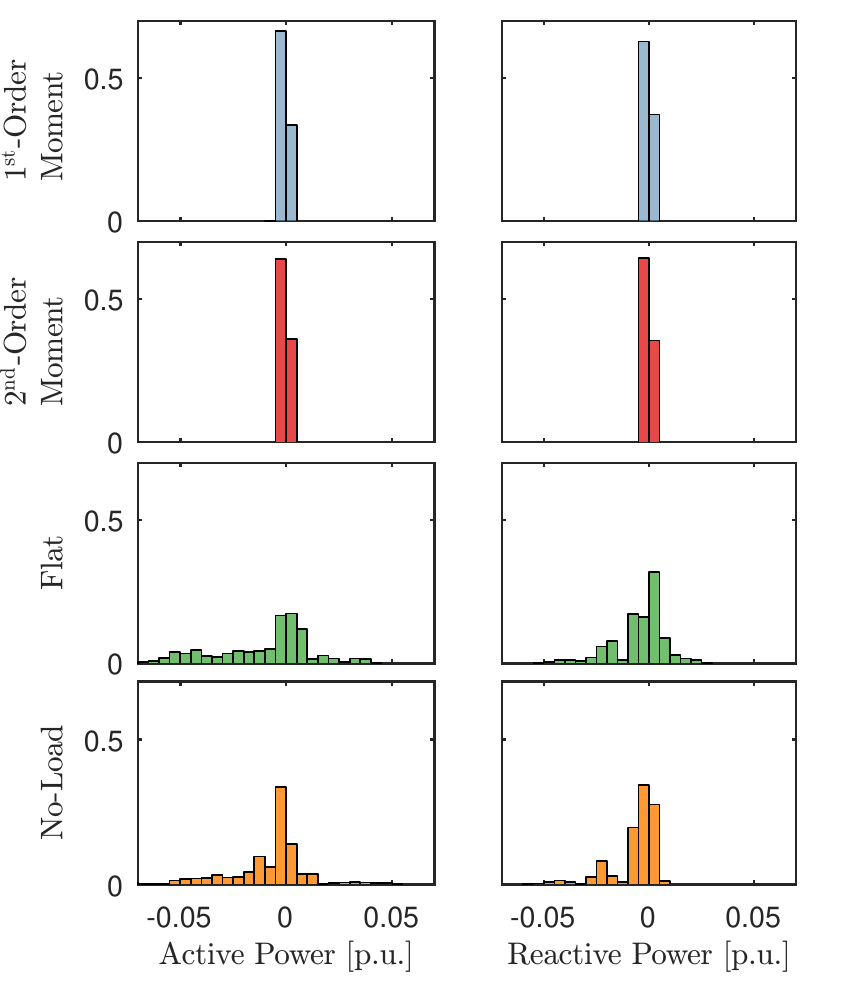}
	\caption{IEEE 33-Bus: Histograms of signed active and reactive power balance errors $p_i$ and $q_i$ for the four linearized models over all buses and 1000 scenarios.}\label{fig:case33}
\end{figure}
\begin{figure}[!t]
	\centering
	\includegraphics[scale=1]{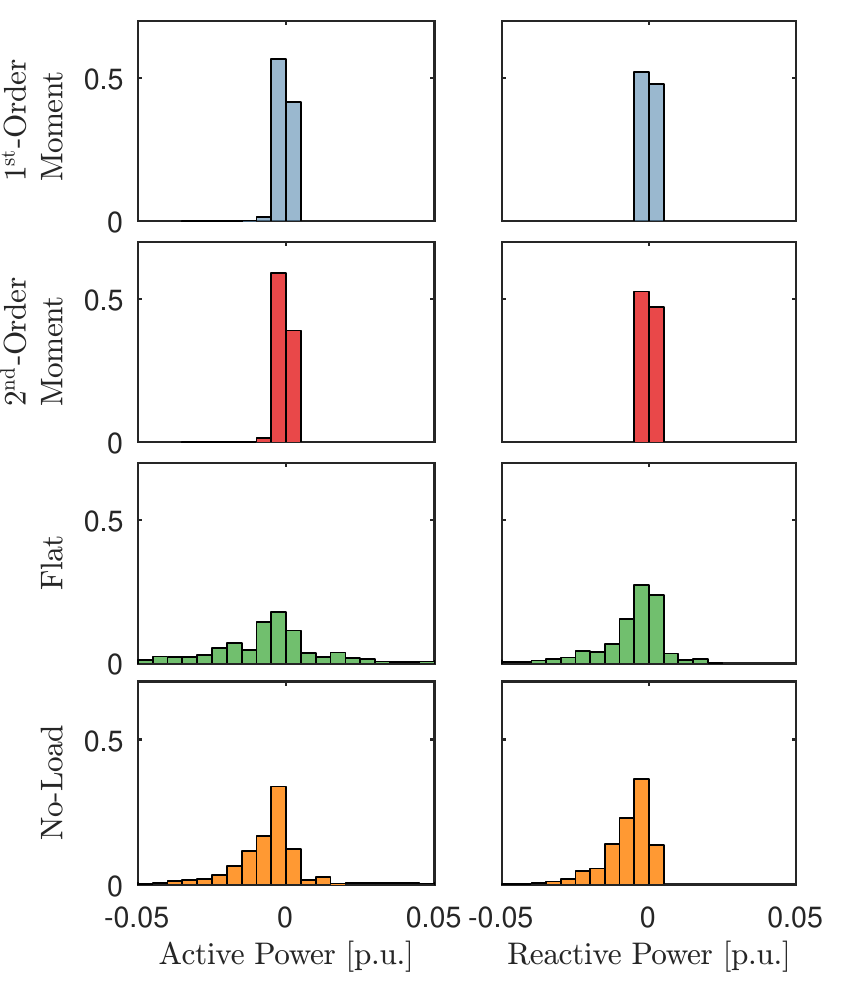}
	\caption{IEEE 57-Bus: Histograms of signed active and reactive power balance errors $p_i$ and $q_i$ for the four linearized models over all buses and 1000 scenarios.}\label{fig:case57}
\end{figure}
\section{Numerical Results}\label{numericals}
We demonstrate the effectiveness of our linearization approach with benchmark networks from the MATPOWER library \cite{Zimmerman2011}, as listed in Table \ref{table:results}. As in \cite{Molzahn2015}, we impose additional apparent power flow limits on all lines so that standard 1st-Order SDP relaxations fail to return an exact solution. The approach is implemented using MATLAB\textsuperscript{TM} and YALMIP \cite{Lofberg2004}. The SDP relaxations and convex approximations are solved with MOSEK\textsuperscript{TM}. The sparsity pattern is detected using \emph{SparsePOP} \cite{Waki2008}.

The following linearization methods are applied to obtain a convex approximation (\ref{eq:convexparam}) and compared in terms of constraint violations and objective value:
\begin{itemize}
	\item The \emph{1st-Order Moment} and \emph{2nd-Order Moment} linearization profiles developed in Section \ref{approxoptimpoint}: Constraints (\ref{eq:convexparamb}) are linearized around the expected value computed by the sparse SDP relaxations (\ref{eq:sparsesdp}) of degree $2k=2$ and $2k=4$.
	\item The \emph{Flat} linearization profile \cite{dhople2015}: Constraints (\ref{eq:convexparamb}) are linearized around the voltage phasor: $\tilde{E}_i=1,\tilde{F}_i=0, \enskip \forall i \in \mathcal{N}$. 
	\item The \emph{No-Load} linearization profile developed in \cite{Bolognani2016} and \cite{dhople2015}, and used in \cite{Guggilam2016} and \cite{DallAnese2017}: Constraints (\ref{eq:convexparamb}) are linearized around $\tilde{E}+j\tilde{F}=-(G+jB)^{-1}(\overline{G}+j\overline{B})(E_0+jF_0)$ where the index $0$ denotes the slack bus. The admittance matrix is decomposed so that the bus current injections can be expressed as:
	\begin{equation*}
	\begin{pmatrix}
	I\\
	I_0
	\end{pmatrix}=\begin{pmatrix}
	G+jB &\overline{G}+j\overline{B}\\
	(\overline{G}+j\overline{B})^T &G_0+jB_0\\
	\end{pmatrix}
	\begin{pmatrix}
	E+jF\\
	E_0+jF_0
	\end{pmatrix}
	\end{equation*}	
\end{itemize}
In current industrial practice, one often assumes that load uncertainty is determined by a small number of latent stochastic variables \cite{taylor2007}. For instance, using the method presented in \cite{Warrington2018}, high-dimensional weather data can be reduced to a low order uncertainty model. In line with this, we consider here the load at each bus as the weighted sum of two latent variables $r_1$ and $r_2$, i.e.~$P_{Li}=P_{Mi}(a_{i,1}r_1+a_{i,2}r_2)$ and $Q_{Li}=Q_{Mi}(a_{i,1}r_1+a_{i,2}r_2)$, where $r=(r_1,r_2)$ is jointly Gaussian $N(\theta,\Theta)$ with $\theta_i=.85$ and $\Theta=
\begin{pmatrix}
0.01&0.002\\
0.002&0.008
\end{pmatrix}$, $P_{Mi}$ and $Q_{Mi}$ are the nominal loads from the MATPOWER library; and $(a_{i,1},a_{i,2})$ are Bus dependent parameters defined as $a_{i,1}=(i-1)/(N_B-1)$ and $a_{i,2}=1-(i-1)/(N_B-1)$. The sampled values $r_1$ and $r_2$ are both truncated to $[0.7,1]$ to maintain feasibility of the linearized models. 

To compute the expected value and standard deviation of the constraint violations and the expected objective value, we solve problem (\ref{eq:convexparam}) with the four different linearization points for $M=1000$ realizations of $r$. We use the same samples to compute the raw moments $z_\gamma$ for the SDP relaxations.

Given an optimal solution $x_{\rm lin}^*(y)$ of a linearized model (\ref{eq:convexparam}) for a load condition $y$, we define the active and reactive absolute constraint violations as $\epsilon_P:=\sum_{i=1}^{N_B} |p_i|$ and $\epsilon_Q:=\sum_{i=1}^{N_B} |q_i|$, where $p_i$ and $q_i$ are the signed active and reactive power balance errors at bus $i$ as defined in (\ref{eq:nodeactive}) and (\ref{eq:nodereactive}). The expected values $\mathbf{E}(\epsilon_P)$ and $\mathbf{E}(\epsilon_Q)$; and standard deviations $\sigma(\epsilon_P)$ and $\sigma(\epsilon_Q)$ are reported in Table \ref{table:results}. No \emph{2nd-Order Moment} linearization is obtained for the 118-bus case because the SDP relaxation with $2k=4$ is not tractable. The \emph{No-Load} profile causes the 118-bus case to become infeasible and is therefore not reported. No inequality constraint violations were observed. 

$\emph{1st-Order Moment}$ and $\emph{2nd-Order Moment}$ significantly outperform $\emph{Flat}$ and $\emph{No-Load}$ in terms of active and reactive constraint violations. In all cases, except for the 9-bus example, \emph{2nd-Order Moment} performs at least as well as \emph{1st-Order Moment} in terms of constraint violations.

We also plot the distribution of active and reactive power balance errors, $p_i$ and $q_i$, over all buses and all $M$ demand scenarios for the 14-bus, 33-bus and 57-bus case in Fig.~\ref{fig:case14}, Fig.~\ref{fig:case33} and Fig.~\ref{fig:case57}. Whereas $\epsilon_P=\sum_i |p_i|$ and $\epsilon_Q=\sum_i |q_i|$ are indicators for the cumulative magnitude of constraint violations, Fig.~\ref{fig:case14}, Fig.~\ref{fig:case33} and Fig.~\ref{fig:case57} show the considerable spread of violations that can occur over all buses for the \emph{Flat} and \emph{No-Load} profiles. The methods $\emph{1st-Order Moment}$ and $\emph{2nd-Order Moment}$ reduce this spread in all cases. Thus, the OPF solutions from the improved linearization lead to fewer real-time operational issues when implemented in reality. For instance, it could lower the need for reserves procurement and deployment, thanks to the reduced power mismatch.

The $\emph{Flat}$ and $\emph{No-Load}$ underestimate the expected cost because their estimates are consistently lower than the lower bounding values provided the SDP relaxations (\ref{eq:sparsesdp}) (see Table~\ref{table:objvalue}). The expected optimal values of $\emph{1st-Order Moment}$ and $\emph{2nd-Order Moment}$ are higher but closer to the lower bounds, suggesting that $\emph{1st-Order Moment}$ and $\emph{2nd-Order Moment}$ provide more accurate cost estimates with lower constraint violations.

The computation times of the dense and sparse SDP relaxations are reported in Table \ref{table:times}. The larger numerical examples can only be solved using the sparse SDP relaxation (\ref{eq:sparsesdp}). The dense SDP relaxation (\ref{eq:primalsdp}) quickly leads to ``out-of-memory'' errors.

\begin{table}[!t]
	\caption{Constraint violations attained over 1000 demand scenarios. Limits on all lines are indicated in parentheses.}\label{table:results}
	\begin{center}
	\begin{tabular}{ l| c c c c c}
	\textbf{Case}&$\mathbf{E}(\epsilon_P)$&$\sigma(\epsilon_P)$&$\mathbf{E}(\epsilon_Q)$ &$\sigma(\epsilon_Q)$&\\
	\hspace{3mm} Method&[p.u.]&[p.u.]&[p.u.]&[p.u.]\\
		\midrule[1pt]
		\input{Table1.tex}
	    \midrule[1pt]
	\end{tabular}
	\end{center}
\end{table}
\begin{table}[!t]
	\caption{Estimates of the Expected Cost in [\textdollar]}\label{table:objvalue}
	\begin{center}
		\begin{tabular}{ l| c c c c c c}
		\multirow{2}{*}{Method}&\multicolumn{6}{c}{Case}\\
		&5-bus&9-bus&14-bus&33-bus&57-bus&118-bus\\
			\midrule[1pt]
			\input{Table3.tex}
			\midrule[1pt]
		\end{tabular}
	\end{center}
\end{table}
\begin{table}[!t]
	\caption{Time taken to solve SDP relaxations to obtain linearization points \emph{1st-Order Moment} and \emph{2nd-Order Moment} using MOSEK\textsuperscript{TM} on an Intel-i5 2.2GHz CPU with 8GB RAM. Cases marked with $\dagger$ are solved on an Intel Xeon E5 2680 2.5GHz and 128 GB RAM.}\label{table:times}
	\begin{center}
		\begin{tabular}{ l| l l l l l l}
			\multirow{2}{*}{SDP Relaxations}&\multicolumn{6}{c}{Case}\\
			&5-bus&9-bus&14-bus&33-bus&57-bus&118-bus\\ 
			\midrule[1pt]
			\textbf{Dense}  &   &   &   &  & & \\
			\hspace{3mm}(\ref{eq:primalsdp}), $2k=2$&0.61s&4.41s&-&-&-&-\\
			\hspace{3mm}(\ref{eq:primalsdp}), $2k=4$&-&-&-&-&-&-\\
			\midrule[1pt]
			\textbf{Sparse}  &   &   &   &  & & \\
			\hspace{3mm}(\ref{eq:sparsesdp}), $2k=2$&0.13s&0.15s&0.69s&0.78s&3.18s&8.83s\\
			\hspace{3mm}(\ref{eq:sparsesdp}), $2k=4$&8.46s&17.84s&390.16s&1.69h$\dagger$&8.97h$\dagger$&-\\
			\midrule[1pt]
		\end{tabular}
	\end{center}
\end{table}
\section{Conclusion and Future work}\label{conclusion}
This paper presented a novel linearization approach that minimizes the constraint violations by taking into account the uncertainty of the demand. The approach relies on the characterization of the optimal operating range of power systems through a hierarchy of sparse SDP relaxations and is shown to outperform other linearization profiles proposed in the literature.

The contribution of this work is limited to static OPFs, but the GMP, on which our linearization approach relies, can be augmented to include dynamic processes, such as the charging of a battery or the temperature propagation in heating and cooling networks, by using occupation measures \cite{Lasserre2007}. The moments of these occupation measures would provide optimal linearization profiles for nonlinear dynamical systems that result from coupling energy networks and storage systems.
\section*{Acknowledgments}
We would like to express our gratitude to Viktor Dorer, Roy Smith and Jan Carmeliet for their valuable help and support. We would also like to thank Georgios Darivianakis, Benjamin Flamm, Mohammad Khosravi and Annika Eichler for fruitful discussions.
\bibliographystyle{IEEEtran}
\bibliography{IEEEabrv,LinearisationPaper}
\appendix
Since the equality constraints of (\ref{eq:opf}) are second degree polynomials, we can restate them using Taylor's theorem: $h_i(x,y)=h_i(x_0,y)+\nabla_x h_i(x_0)(x-x_0)+(x-x_0)^T H_{i} (x-x_0)$, where $H_i$ is a constant matrix representing the Hessian. Selecting the linearization $h_i^{\rm lin}(x,y)$ to match the constant and linear terms of the Taylor approximation, leaving only the quadratic term, and using (\cite{Lasserre2010}, Corollary 2.3), the minimization (\ref{eq:meanconstraint}) can be equivalently written as:
\begin{align}\label{eq:meanconstrainthessian}
\min_{x_0} \int_{\mathbf{K}} \sum_{i=1}^{N_h}|(x-x_0)^T H_{i} (x-x_0)| d\mu^*
\end{align}
We distinguish between two cases: all Hessians are semi-definite or at least one Hessian is indefinite. For the case of positive semi-definite and negative semi-definite Hessians, denoted by $H^+$ and $H^-$, the problem (\ref{eq:meanconstrainthessian}) is a convex optimization for which an exact closed-form solution of can be found. An optimal linearization point $x_0^* \in \mathbb{R}^n$ is such that $\frac{\partial}{\partial x_0} \int_{\mathbf{K}} \sum_s(x-x_0)^T H_{s}^+ (x-x_0)-\sum_t(x-x_0)^T H_{t}^- (x-x_0) d\mu=0$ when $x_0=x_0^*$, where the equality constraints $h_i(x,y)$ are divided into constraints $h_s(x,y)$ with positive semi-definite Hessians  $H_{s}^+$ and constraints $h_t(x,y)$ with negative semi-definite Hessians $H_{t}^-$.

Leibniz's integral rule allows us to take the derivative of (\ref{eq:meanconstrainthessian}) with respect to $x_0$:
\begin{equation}\label{eq:meanconstraintderivation}
\begin{aligned}
\frac{\partial}{\partial x_0} \int_{\mathbf{K}}& \sum_s(x-x_0)^T H_{s}^+ (x-x_0)\\&-\sum_t(x-x_0)^T H_{t}^- (x-x_0) d\mu^*\\
=\int_{\mathbf{K}} &\frac{\partial}{\partial x_0} \bigg(\sum_s(x-x_0)^T H_{s}^+ (x-x_0)\\
&-\sum_t(x-x_0)^T H_{t}^- (x-x_0)\bigg)d\mu^*\\
=\int_{\mathbf{K}} &\sum_s 2 H_{s}^+(x-x_0)-\sum_t 2 H_{t}^-(x-x_0)d\mu^*\\
=\Big(&\sum_s 2 H_{s}^+-\sum_t 2 H_{t}^-\Big)\Big(\mathbf{E}_{\mu^*}(x)-x_0\Big) = \mathbf{0},    
\end{aligned}
\end{equation}
where $\mathbf{0}$ is a zero vector. Thus, the expected absolute constraint violation is minimized when $x_0^*=\mathbf{E}_{\mu^*}(x)$ for the case of positive and negative semi-definite Hessians.

For the case that at least one of the Hessians is indefinite, a convex approximation of (\ref{eq:meanconstrainthessian}) is established. Any indefinite matrix can be decomposed into a difference of positive semi-definite matrices: $H_{i}=\tilde{H}^+_{i}-\tilde{H}^-_{i}$ where $\tilde{H}^+_i\succeq 0$ and $\tilde{H}^-_i\succeq 0$.
If one decomposes the indefinite Hessians in (\ref{eq:meanconstrainthessian}) and applies the triangle inequality, a convex upper approximation can be derived:
\begin{equation}\label{eq:meanconstrainthessianindef}
\begin{aligned}
\min_{x_0} \int_{\mathbf{K}} &\sum_i|(x-x_0)^T \tilde{H}^+_{i}(x-x_0)\\
&-(x-x_0)^T \tilde{H}^-_{i} (x-x_0)| d\mu^*\\
\leq\min_{x_0} &\int_{\mathbf{K}} \sum_i (x-x_0)^T \tilde{H}^+_{i}(x-x_0)\\
&+(x-x_0)^T \tilde{H}^-_{i} (x-x_0)d\mu^*
\end{aligned}
\end{equation}
Analogous to (\ref{eq:meanconstraintderivation}), taking the derivative with respect to $x_0$ of the right-hand side of (\ref{eq:meanconstrainthessianindef}) leads to the closed-form expression $x_0^*=\mathbf{E}_{\mu^*}(x)$.

\ifCLASSOPTIONcaptionsoff
  \newpage
\fi

\end{document}

%% file: Table1.tex
\textbf{5-Bus} &   &   &   &   \\
\hspace{3mm} 1st -Order Moment & 0.008 & 0.015 & 0.023 & 0.043 \\
\hspace{3mm} 2nd-Order Moment & 0.008 & 0.014 & 0.022 & 0.041 \\
\hspace{3mm} Flat & 0.359 & 0.028 & 0.787 & 0.018 \\
\hspace{3mm} No-Load & 0.361 & 0.030 & 0.698 & 0.028 \\
\midrule
\textbf{9-Bus} (120 VA) &   &   &   &   \\
\hspace{3mm} 1st -Order Moment & 0.004 & 0.005 & 0.002 & 0.002 \\
\hspace{3mm} 2nd-Order Moment & 0.006 & 0.003 & 0.003 & 0.002 \\
\hspace{3mm} Flat & 0.277 & 0.012 & 0.241 & 0.023 \\
\hspace{3mm} No-Load & 0.811 & 0.027 & 0.609 & 0.058 \\
\midrule
\textbf{IEEE 14-Bus} (25 VA) &   &   &   &   \\
\hspace{3mm} 1st -Order Moment & 0.004 & 0.001 & 0.005 & 0.002 \\
\hspace{3mm} 2nd-Order Moment & 0.001 & 0.002 & 0.002 & 0.002 \\
\hspace{3mm} Flat & 0.114 & 0.004 & 0.113 & 0.006 \\
\hspace{3mm} No-Load & 0.074 & 0.006 & 0.144 & 0.016 \\
\midrule
\textbf{IEEE 30-Bus} (130 VA) &   &   &   &   \\
\hspace{3mm} 1st -Order Moment & 0.002 & 0.002 & 0.001 & 0.001 \\
\hspace{3mm} 2nd-Order Moment & 0.002 & 0.002 & 0.001 & 0.001 \\
\hspace{3mm} Flat & 0.668 & 0.091 & 0.253 & 0.033 \\
\hspace{3mm} No-Load & 0.551 & 0.082 & 0.198 & 0.025 \\
\midrule
\textbf{IEEE 57-Bus} (77 VA) &   &   &   &   \\
\hspace{3mm} 1st -Order Moment & 0.020 & 0.025 & 0.004 & 0.005 \\
\hspace{3mm} 2nd-Order Moment & 0.021 & 0.023 & 0.004 & 0.005 \\
\hspace{3mm} Flat & 1.289 & 0.170 & 0.745 & 0.106 \\
\hspace{3mm} No-Load & 0.746 & 0.034 & 0.401 & 0.043 \\
\midrule
\textbf{IEEE 118-Bus} (110 VA) &   &   &   &   \\
\hspace{3mm} 1st -Order Moment & 0.467 & 0.133 & 0.344 & 0.099 \\
\hspace{3mm} 2nd-Order Moment & \multicolumn{4}{c}{intractable}\\
\hspace{3mm} Flat & 6.832 & 1.820 & 4.851 & 1.326 \\
\hspace{3mm} No-Load & \multicolumn{4}{c}{infeasible}\\

%% file: Table3.tex
\textbf{SDP Relaxations}  &   &   &   &  & & \\
\hspace{3mm}Sparse (\ref{eq:sparsesdp}), $2k=2$& 10532 & 4214 & 7673 & 7236 & 34400 & 108410 \\
\hspace{3mm}Sparse (\ref{eq:sparsesdp}), $2k=4$& 10533 & 4227 & 7675 & 7238 & 34453 & - \\
\midrule
\textbf{Linearizations}  &   &   &   &  & & \\
\hspace{3mm}1st -Order Moment & 10667 & 4240 & 7677 & 7253 & 34509 & 109395 \\
\hspace{3mm}2nd-Order Moment & 10667 & 4240 & 7677 & 7253 & 34510 & - \\
\hspace{3mm}Flat & 10459 & 4179 & 7612 & 6763 & 34086 & 100825 \\
\hspace{3mm}No-Load & 10460 & 4097 & 7604 & 6738 & 33960 & - \\